\documentclass{elsart}

\usepackage{graphicx}
\usepackage{amssymb}

\newcommand {\BM} {basilar membrane}
\newcommand {\micron} {$\mu$m}
\newcommand {\Bekesy} {B\'{e}k\'{e}sy}

\begin{document}

\begin{frontmatter}



\title{
A Comprehensive Three-Dimensional Model of the Cochlea
\thanksref{grant_ref}
}

\thanks[grant_ref]
{Part of this research was supported by NIH grant RO1-DC04084.}


\author[umich]{Edward Givelberg}

\address[umich]{
Department of Mechanical Engineering,
University of Michigan, Ann Arbor
}

\author[caltech]{Julian Bunn}

\address[caltech]{
Center for Advanced Computing Research,
California Institute of Technology, Pasadena, USA
}

\begin{abstract}
The human cochlea is a remarkable device, able to discern extremely
small amplitude sound pressure waves, and discriminate between very
close frequencies. 
Simulation of the cochlea is computationally challenging due
to its complex geometry, intricate construction and small physical
size. 
We have developed, and are continuing to refine, a detailed 
three-dimensional computational model based on an accurate
cochlear geometry obtained from physical measurements. In the model,
the immersed boundary method is used to calculate the fluid-structure 
interactions produced in response to incoming sound waves. The model
includes a detailed and realistic description of the various elastic 
structures present.

In this paper, we describe the computational model and its
performance on the latest generation of shared memory 
servers from Hewlett Packard. Using compiler
generated threads and OpenMP directives, we have achieved
a high degree of parallelism in the executable, which has 
made possible several large scale numerical simulation experiments 
that study the interesting features of the cochlear system. 
We show several results from these simulations, 
reproducing some of the basic known characteristics of cochlear
mechanics.
% with the currently accepted theories of how the system works.
% We show several results from these simulations, and contrast them
% with the currently accepted theories of how the system works.

\end{abstract}

\end{frontmatter}

%============================================================================
\section{Introduction}
\label{sec:introduction}
The cochlea is the part of the inner ear where acoustic signals are
transformed into neural pulses which are then signaled to the brain.
It is a small organ, the size of a child's marble, which is a miracle
of nature's engineering.
It is sensitive to signals ranging over more than seven orders 
of magnitude, from a whisper to an explosion.
% The human cochlea 
We can hear sounds ranging in frequency from approximately 10 Hz to
about 20 kHz, 
and distinguish between tones whose frequencies differ 
by less than a percent.
Trained musicians, for example, are capable of differentiating
between pure tones of 1000 Hz and 1001 Hz.

These remarkable signal processing capabilities are achieved by
a complicated mechanism involving 
both interaction of elastic material components immersed in fluid
% fluid-structure interactions 
and electro-chemical processes.
After decades of research a fascinating picture of the 
cochlear mechanics has emerged, but
% Despite decades of research and much progress, 
the precise nature of the mechanisms
responsible for the extreme sensitivity, sharp frequency selectivity
and the wide dynamic range of the cochlea remains unknown.
The goal of our project is to build computational modeling tools 
that we hope will assist in the understanding of the cochlea works.
We have constructed a comprehensive three-dimensional computational
model of the fluid-structure interactions of the cochlea using the
immersed boundary method.
This is the first model that incorporates 
the curved cochlear anatomy 
based on real physical measurements,
uses the non-linear Navier-Stokes equations 
of viscous fluid dynamics
and 
includes a detailed and realistic description of the various elastic 
structures present.
% The model we have constructed 
% Our model is the only model that incorporates
% the curved three-dimensional geometry of the cochlea,
% uses the viscous non-linear Navier-Stokes fluid equations
% and incorporates detailed models of the elastic properties
% of various cochlear components.
For example, the helicoidal
basilar membrane is modeled as an elastic shell
using partial differential equations.
% uses a shell theory to model the basilar membrane
% and incorporates detailed elastic models of the other materials.

We have developed a suite of software codes to support 
our studies of the cochlea using the immersed boundary method. 
These include codes for the generation of simulation input models 
(implemented in C++),
the main immersed boundary numerical solver,
and those for the analysis and visualization of results 
(implemented in Java and C++).
The main workhorse in this suite is the general purpose
immersed boundary code, which is written in C and requires extensive 
computing resources 
(CPU power, memory, and available disk space for storage of the
results files). 

Our present work builds on the first author's implementation of the
immersed boundary method for elastic shells \cite{G}.
In this paper we outline the construction of the cochlea model
and present preliminary results of several large scale numerical
experiments.
These experiments reproduce some of the basic features 
of cochlear mechanics and demonstrate the promise of large scale 
computational modeling approach to answering
important questions about nature's miraculous engineering
design of the cochlea.

The rest of the paper is organized as follows.
Section \ref{sec:cochlea} is a short introduction
% we mention some of the basic facts of
to cochlear mechanics,
% We also briefly survey 
including a brief survey of some of the previous work on
cochlear modeling.
For a comprehensive up to date introduction to the subject
see Geisler's excellent book \cite{Geisler1998}.
Section \ref{sec:ibm} describes the general methodology of the
immersed boundary method.
We write down the fluid-structure interaction equations
and outline the general first-order numerical method to solve 
these equations.
The method has been extensively tested for different models of
the elastic immersed boundary.
For more detailed presentations see \cite{Peskin1994}, \cite{Peskin2002}
and \cite{G}.
We also discuss briefly the implementation of the immersed boundary
algorithm and its optimization.
Section \ref{sec:model} describes the construction of the cochlea model.
% (anatomy, software, validation, visualization).
In section \ref{sec:results} we present results of several numerical
simulations.
Our cochlea model is a work in progress and
we conclude with an outline of future directions for this project.

\section{The Cochlea}
\label{sec:cochlea}
\subsection{Cochlear Mechanics}
The cochlea is a small snail-shell-like cavity in the temporal bone,
which has two openings, the oval window and the round window.
The cavity is filled with fluid and is sealed by two elastic membranes that
cover the windows.
The spiral canal of the cochlea
is divided lengthwise into two passages 
by the cochlear partition consisting of the bony shelf and
the so-called basilar membrane.
These passages are the scala vestibuli and scala tympani,
and they
connect with each other at the apex, called the helicotrema.
External sounds set the ear drum in motion,
which is conveyed to the inner ear by the ossicles,
three small bones of the middle ear, the malleus, incus and stapes.
The ossicles function as an impedance matching device,
focusing the energy of the ear drum on the oval window of the
cochlea.
The piston-like motion of the stapes against the oval window
displaces the fluid of the cochlea, so generating traveling waves
that propagate along the basilar membrane.

% Practically everything 
Much of what
we know about the waves in the cochlea 
was discovered in the 1940s 
by Georg von \Bekesy\ \cite{vonBekesy},
% Von \Bekesy\ 
who carried out experiments in cochleae extracted 
from human cadavers.
Von \Bekesy\ studied the cochlea as a passive mechanical filter
that utilizes a system of elastic components immersed in a fluid
for analysis of incoming sounds.
He observed that
a pure tone input generates a 
traveling wave propagating along the basilar membrane.
The wave amplitude rises gradually,
reaching a peak at a specific location along the membrane,
after which it decays rapidly.
The peak location depends on the frequency of the input tone,
with high frequencies peaking close to the stapes, 
and the lower frequencies further towards the apex.
% and the lower the input frequency, the further the location
% of its peak from the stapes.
This ``place principle'' is a crucial mechanism
of frequency analysis in the cochlea.
Resting on the basilar membrane is the microscopic organ of Corti,
a complicated structure containing sensory receptors 
called hair cells.
The hair cells detect fluid motion and 
provide input to the afferent nerve fibers
that send action potentials
to the brain.
% The macro-mechanical fluid-shell system performs frequency analysis
% since 
Thus, a pure tone input activates only 
a specific group of
% a narrow band of
hair cells depending on its frequency,
with the 
characteristic frequency locations monotonically decreasing 
along the basilar membrane
from the stapes to the helicotrema.

Von \Bekesy\ found that the
%'s experiments revealed the important role the
\BM's elastic properties play an important role
in the wave mechanics of the cochlea.
% with the cochlea and the various test
% models he has constructed revealed
% Von \Bekesy\ observed that
Despite its name, the basilar membrane 
is in fact an elastic shell whose
compliance increases exponentially with length.
(Unlike an elastic membrane, an elastic shell is not under inner
tension, i.e. when it is cut the edges do not pull apart.
The compliance of an elastic shell is defined as the amount
of volume displaced per unit pressure difference 
across the surface of the shell.)
% The \BM's elastic properties 
% play an important role in the wave mechanics of the cochlea.
Von \Bekesy's experiments further revealed that 
the wave 
propagates in the basilar membrane in the direction of increasing
compliance regardless of the location of the source of excitation in
the fluid.
This phenomenon came to be known as ``the traveling wave paradox''
and it is very important because a significant part 
of our hearing depends on
bone conduction, where the input to the cochlea is coming not through
the stapes, but through the vibration of the bony walls.
Bone conduction is easily demonstrated by placing 
a vibrating tuning fork in contact with the forehead,
resulting in the subject hearing the tone of the fork frequency.
%
% %
% \begin{figure}[h]
% % \setlength{\epsfxsize}{5in}
% % \hspace{0.5in}
% \vspace{0.2in}
% % \epsffile{corti.ps}
% \includegraphics{corti.pdf}
% \caption{A schematic diagram of a cross-section of the organ of Corti.}
% \label{fig:Corticross}
% \end{figure}
% %

Mathematically,
the macro-mechanical system of the cochlea
can be described by the Navier-Stokes
equations of incompressible fluid mechanics coupled with equations
modeling the elastic properties 
of the basilar membrane and the membranes 
of the oval and the round windows.
% In fact, measurements of these properties seem to vary by an
% order of magnitude [][][][],
% but assume that they are known.
The mathematical problem of solving this system of
partial differential equations on a three-dimensional domain 
with intricate curved geometry
is very difficult.
Since the displacements of the basilar membrane are extremely small
(on a nanometer scale), 
the system operates in a linear regime.
Analysis shows that the waves in the cochlea
resemble shallow water waves
% , like the ripples on the surface of a pond
\cite{Leveque1988}.
% The passive macro-mechanical filter can be shown 
% approximately to perform a
% wavelet transform.

%
While the macro-mechanics of the cochlea break up the incoming
sound into its frequency components,
it was suggested
as early as 1948 that 
a passive mechanical system alone
cannot explain the extreme sensitivity and 
frequency selectivity of the cochlea;
some kind of amplification is necessary \cite{Gold1948}.
Indeed, analysis of cochlear macro-mechanics indicates that the
traveling wave focusing is not sufficiently sharp, with 
some estimates suggesting that, at its threshold, human hearing
is about a hundred times more sensitive than 
what would be expected
from a passive macro-mechanical filter of the cochlea.

In 1967 Johnstone and Boyle \cite{Johnstone1967} 
utilized the M\"{o}ssbauer effect
to carry out measurements of the basilar membrane vibrations 
{\em in vivo},
with more than 100-fold improvement in resolution over von \Bekesy's
measurements.
These measurements revealed that in the live cochlea the peak of the
wave envelope is, in fact, well localized \cite{Rhode1971}.
Since then
the fine tuning of the cochlea has been confirmed experimentally
in many studies,
e.g. \cite {Robles1984,Sellick1982,Khanna1982}.
%
% a pure tone input activates only a narrow band of hair cells.
% There is evidence that humans can detect auditory stimuli that provide
% each hair cell with an energy near the thermal level [DeVries].
% The ear can respond to hair-bundle stimuli less than 1 nm 
% in peak-to-peak magnitude [Sellick, Patuzzi...]
% and can discriminate frequencies that differ by a fraction of 1\%.
This led to attempts to understand the mechanical 
and electro-chemical processes
at the level of the hair cells and to a renewed interest in
the conjecture that the cochlea
contains an energy source and acts as a mechanical amplifier.

Further indirect support for the active amplification hypothesis
was found
with the discovery of the existence of otoacoustic emissions
% The SOAEs were discovered in 
\cite {Kohlloffel1972,Kemp1978}.
Otoacoustic emissions may be recorded during or after acoustical
stimulation using a sensitive microphone placed close to the ear drum.
Emissions can also be detected when electric current is applied
to the cochlea
\cite{Nuttall1995}.
Spontaneous otoacoustic emissions (SOAEs) occur in humans and other
species.
% They have been found to occur in more than a half of the human
% population and 
In the severe cases
SOAEs are the cause of {\it objective tinnitus},
a common complaint of patients with ``ringing ears''.
% Recent studies have found that otoacoustic emissions occur not only
% in people suffering from tinnitus, but in more than half of the 
% general population.
Experiments show that a SOAE can be suppressed when a stimulus tone 
is presented at a nearby frequency.
An isosuppression tuning curve is a curve obtained
by measuring the otoacoustic emission while
varying the amplitude and the frequency of the stimulus.
In both mammalian and nonmammalian cochleae SOAEs 
and the process responsible for sharp frequency selectivity
display similar characteristics.
Measurements show that the isosuppression tuning curve
closely resembles an ordinary tuning curve 
\cite{Clark1984},
which measures responsiveness to acoustical stimulus,
leading to the conjecture that the cochlear amplifier is also
the source of the SOAEs.

% Complex sounds composed of several pure tones typically evoke
% a basilar membrane response that is similar to a superposition of
% the membrane's responses to the constituent frequencies,
% however, in some experiments the basilar membrane response 
% was shown to be more complicated.
% (For example, in response to certain combinations of two pure tones
% the basilar membrane generates the so-called distortion product,
% causing the cochlea to respond as though three pure tone input
% frequencies are present.)

Presently,
understanding the nature of the active mechanism in the cochlea
is the main open problem of hearing research.
The live cochlea is a remarkable highly non-linear filter,
but its function crucially depends on the underlying linear filter
of the passive cochlear mechanics,
which is still not sufficiently well understood.
Basic questions about the role of geometry and the elastic properties of
the basilar membrane in cochlear mechanics remain open.
Answering such questions is important not only to our understanding of
the cochlea, but also in solving important engineering problems,
such as a cochlear transducer design (see \cite{Ramamoorthy2002}).

\subsection{Cochlear Models}

% Direct measurements of the basilar membrane in vivo 
% are very difficult.
% A computational model of the cochlea is needed because
% make it very difficult to measure the vibrations of the basilar
% membrane in vivo.
% Computational modeling of the cochlea is therefore necessary to aid
% in cochlear research.

% Much of the difficulty in cochlear research is due to
% the small size and the complex geometry of the cochlea.
Extensive research in cochlear modeling has been carried out
over the years.
Because of mathematical difficulty mostly 
simplified one or two-dimensional models
that sought to incorporate
some 
% of cochlear features
aspects of cochlear mechanics
have been proposed.

Early one-dimensional ``transmission line'' model of the cochlea
\cite {Zwislocki1948}, 
\cite {Peterson1950}, 
\cite {Fletcher1951}
has assumed that the fluid pressure is constant over a cross section
of the cochlear channel.
% which leads to an ordinary differential
% equation for the fluid pressure. 
The fluid is assumed to be incompressible and inviscid,
and the
basilar membrane is modelled as a damped, forced harmonic
oscillator with no elastic coupling along its length.
Qualitatively, this model has been shown to capture the basic features
of the basilar membrane response.
Quantitatively, however, it yields large discrepancies with
measurement results \cite {Zweig1976}.

Two-dimensional models by Ranke \cite {Ranke1950}
and Zwislocki \cite {Zwislocki1965} 
make similar assumptions on the cochlear fluid and 
the basilar membrane.
Ranke's model uses a deep water approximation,
while Zwislocki used the
shallow water theory in his model.
These models were further developed in
% by many authors 
\cite {Siebert1974}, 
\cite {Lesser1972},
\cite {Allen1979},
\cite {Allen1977}
and other works.
The most rigorous analysis
of a two-dimensional model with fluid viscosity 
was carried out by
Leveque, Peskin and Lax in \cite{Leveque1988}.
Their cochlea 
% A theoretical analysis of a two dimensional flat model 
% of the cochlea
% was carried out by Leveque, Peskin and Lax in \cite{LPL}.
% In this model the cochlea is represented by 
is represented by
a plane (i.e. a strip of infinite length and infinite depth)
and the basilar membrane, by an infinite line dividing the plane into
two halves.
% An asymptotic analysis shows that this
% model captures the main 
% features of the cochlear wave mechanics that were found by von Bekesy.
The linearized equations are reduced to a functional equation by
applying the Fourier transform in the direction parallel to the
basilar membrane and then solving the resulting ordinary differential
equations in the normal direction.  The functional equation derived in
this way is solved analytically, and the solution is evaluated both 
numerically and also asymptotically (by the method of stationary phase).
This analysis reveals that the waves in the cochlea resemble
shallow water waves, i.e. ripples on the surface of a pond.
A distinctive feature
of this paper is the (then speculative) consideration 
of negative basilar 
membrane friction, i.e., of an amplification mechanism operating within
the cochlea.

Other two-dimensional models incorporate more sophisticated
representations of the \BM,
using, for example,  elastic beam and plate theory
(\cite {Bogert1951},
\cite {Schroeder1973},
\cite {Kim1973},
\cite {Allaire1974},
\cite {Steele1974},
\cite {Inselberg1976},
\cite {Chadwick1980},
\cite {Holmes1982}).
Three-dimensional models were considered by Steel and Taber 
\cite {Steele1979} and de Boer \cite {deBoer1982}, who used asymptotic
methods and computations,
obtaining a slightly improved fit of the experimental data.
Their experience seems to indicate that geometry may play 
a significant role in the problem.
In particular, the effect of the spiral coiling of the cochlea on the
wave dynamics remains unresolved. 
It has been considered by several authors
(see
\cite {Viergever1978},
\cite {Loh1983},
\cite {Steele1985},
\cite {Manoussaki2000}).
% The prevailing opinion is that the coiling has only a marginal
% effect on the \BM\
% motion.

With the development of more powerful computers it became possible
to construct more detailed computational models of the cochlea.
A two-dimensional computational model of the cochlea 
was constructed by Beyer \cite{Beyer1992}.
In this model the cochlea is a flat rectangular strip divided 
into two equal halves by a line which represents the basilar membrane.
The fluid is
modelled by the full Navier-Stokes equations with viscosity terms, 
but elastic coupling along the \BM\
is not incorporated.
Beyer has used a modification of Peskin's immersed boundary method,
originally developed for modeling the fluid dynamics of the heart
\cite{Peskin1977}.
Several
three-dimensional computational models have been reported,
such as Kolston's model
\cite{Kolston1996,Kolston1999},
intended to simulate
the micro-mechanics of the cochlear partition
in the linear regime (i.e. near the threshold of hearing),
and Parthasarati, Grosh and Nuttal's
\cite{Parthasarati2000}
hybrid analytical-computational model
using WKB approximations and finite-element methods.

% The model we present in this paper
% is the first comprehensive three-dimensional
% computational model of 
% the passive cochlear macro-mechanics.
% In the future 
% % we hope to incorporate features
% we will extend this model to 
% incorporate features of the organ of Corti
% in order to study the active mechanism of the cochlea.

% The study of the active mechanism of the cochlea is the ultimate goal 
% of our project.
% This paper describes the preliminary step towards that goal, namely 
% the construction of a comprehensive three-dimensional 
% computational model of 
% the passive cochlear macro-mechanics.
%
% Cochlear mechanics is a vast subject of study with thousands of 
% experimental and modeling works reported in the literature.
% Because of the mathematical complexity previous modeling efforts have
% largely treated one or two-dimensional system using asymptotic
% or hybrid asymptotic-computational models.
%
% The model we have constructed is the only model that incorporates
% the full curved three-dimensional geometry of the cochlea,
% uses the viscous non-linear Navier-Stokes fluid equations,
% uses a shell theory to model the basilar membrane
% and incorporates detailed elastic models of the other materials.

\section{The Immersed Boundary Method}
\label{sec:ibm}

The primary method used in our construction of the cochlea model is the
immersed boundary method of Peskin and McQueen.
% \cite{Peskin1994}.
% The immersed boundary method of Peskin and McQueen \cite{Peskin1994} 
It is a general
numerical method for modeling an elastic material
immersed in a viscous incompressible fluid. 
It has proved to be particularly useful for
computer simulation of various biofluid dynamic systems.
In this section we outline the general framework of the immersed boundary method.
So far most applications of the method have treated
the elastic (and possibly active) biological tissue
as a collection of elastic fibers 
immersed in a viscous incompressible fluid.
For details of
this formulation of the method together with references 
to many applications see \cite{Peskin1994} and \cite{Peskin2002}.
% The most important component of the cochlea is the
% so-called ``basilar membrane''.
%
The immersed boundary framework is suitable
for modeling not only elastic fibers, 
but also different elastic materials having
complicated structure. 
The cochlea model makes an essential use of the
immersed boundary method for shells,
as developed in
% incorporates a shell theory into the immersed boundary framework,
% obtaining a practical computational method for modeling
% an elastic shell immersed in fluid 
\cite{G}.
The main advantage of 
the immersed boundary method is its conceptual simplicity:
the viscous incompressible fluid is described 
by the Navier-Stokes equations, 
the geometry of the model 
mirrors the real-life curved three-dimensional
cochlear geometry and models for elastic and active material
components can be naturally integrated.

\subsection {The Equations of the Model}

\def\Rthree{{\bf R}^3}
\def\Rtwo{{\bf R}^2}
\def\boldu{{\bf u}}
\def\boldF{{\bf F}}
\def\boldf{{\bf f}}
\def\boldx{{\bf x}}
\def\boldX{{\bf X}}
\def\boldq{{\bf q}}
\def\dubydt{\frac {\partial \boldu} {\partial t}}
\def\dXbydt{\frac {\partial \boldX} {\partial t}}

The immersed boundary method is based on a Lagrangean  formulation 
of the fluid-immersed material system.
The fluid is described in the standard cartesian coordinates
on $\Rthree$,
while
the immersed material is described in a different curvilinear 
coordinate system.
%
% We first turn to the Navier-Stokes equations.
Let $\rho$ and $\mu$ denote the density and the viscosity of
the fluid,
and let ${\boldu}({\boldx}, t)$ 
and $p({\boldx}, t)$ denote its velocity and pressure,
respectively.
The Navier-Stokes equations of
a viscous incompressible fluid are:
\begin{eqnarray}
\label {eq:NS1}
\rho \left( \dubydt + \boldu \cdot \nabla \boldu \right)
        & = & \mbox{} - \nabla p + \mu \nabla^2 \boldu + \boldF
\\
\label {eq:NS2}
\nabla \cdot \boldu & = & 0,
\end{eqnarray}
where $\boldF$ denotes the density of the body force
acting on the fluid.
% the immersed material applies on the fluid.
For example, if the immersed material
is modeled as a thin shell,
then $\boldF$ is a singular vector field,
which is zero everywhere, except possibly on the surface
representing the shell.
The numerical method uses a discretization of the
Navier-Stokes equations (\ref{eq:NS1}) and (\ref{eq:NS2}) on
a periodic rectangular grid.

Let $\boldX(\boldq, t)$ denote the position of the immersed material
in $\Rthree$.
% In the case of
For a shell, $\boldq$ takes values in a domain
$\Omega \subset \Rtwo$,
and $\boldX(\boldq, t)$ is a 1-parameter family of
surfaces indexed by $t$,
i.e., $\boldX(\boldq, t)$ is the middle surface
of the shell at time $t$.
Let $\boldf(\boldq, t)$ denote the force density that
the immersed material applies on the fluid.
% Then Newton's first law can be expressed as
Then
\begin{eqnarray}
\label {eq:fluidF}
\boldF(\boldx, t) = \int \boldf(\boldq, t)\delta(\boldx
- \boldX(\boldq, t)) \,d\boldq,
\end{eqnarray}
where $\delta$ is the Dirac delta function on $\Rthree$.
This equation merely says that the fluid feels the force that the
immersed material exerts on it, 
but it is important in the numerical method,
where it is one of the equations determining fluid-material
interaction.
The other interaction equation is the
no-slip condition for a viscous fluid:
\begin{eqnarray}
\label {eq:fluidU}
\dXbydt & = & \boldu (\boldX (\boldq, t), t)
\nonumber \\
& = & \int \boldu(\boldx, t) \delta(\boldx - \boldX(\boldq, t))
\, d\boldx.
\end{eqnarray}
\def\boldN{{\bf N}}
\def\boldW{{\bf W}}
The system has to be completed by specifying
the force $\boldf(\boldq, t)$ of the immersed material.
In a complicated system such as the cochlea
the immersed material consists of many different components:
membranes, bony walls, an elastic shell representing the
basilar membrane,
and various cells of the organ of Corti,
including outer hair cells, which may actively generate forces.
For each such component it is  necessary to specify its own
computation grid and an algorithm to compute its force $\boldf$.
It is in the specification of these forces that
models for various system components integrate into
the macro-mechanical model.

% \subsection*{The immersed boundary method: implementation.}
\subsection{The Numerical Method}
We describe here 
a first-order immersed boundary numerical scheme,
which
is the easiest
to implement and has the important advantage of modularity:
incorporating various models of
immersed elastic material is straightforward.
% in the current implementation
% solving the fully non-linear Navier-Stokes equations
% is not more expensive than solving linearized equations.
%
\def\deltat{{\Delta t}}

Let $\deltat$ denote the duration of a time step.
It will be convenient to denote the time step by the superscript.
For example
$ \boldu^n(\boldx) = \boldu(\boldx, n \deltat). $
At the beginning of the $n$-th time step $\boldX^n$ and $\boldu^n$
are known.
Each time step proceeds as follows. 
\begin {enumerate}
\item
Compute the force $\boldf^n$ that the immersed boundary
applies to the fluid.
For simple materials, such as fibers, 
this is a straightforward computation (see \cite{Peskin1994}).
% For an elastic shell this computation is more complicated.
For a detailed description of a shell immersed boundary force
computation see \cite{G}.
\item
% Use the interpolation equations 
Use (\ref {eq:fluidF})
to compute the external force on the 
fluid $\boldF^n$.
\item
Compute the new fluid velocity $\boldu^{n + 1}$ from the Navier Stokes
equations.
\item
% Using the interpolation equations 
Use (\ref {eq:fluidU})
to compute the new position 
$\boldX^{n + 1}$ of the immersed material.
\end {enumerate}
%
% The computation of the force in step 1 depends on the material being modeled.
% For the case of an elastic shell please see \cite{G}.
We shall now describe in detail the computations in steps 2 --- 4,
beginning with the Navier-Stokes equations.

The fluid equations are discretized on a rectangular lattice of mesh
width $h$. 
We will make use of the following difference operators which act on
functions defined on this lattice:
\def\Dpi{{D^+_i}}
\def\Dmi{{D^-_i}}
\def\Dci{{D^0_i}}
\def\boldD0{{\bf D^0}}
\def\bolde{{\bf e}}
\begin{eqnarray}
\label {eq:differenceOpsFirst}
\Dpi \phi (\boldx) & = &
	\frac {\phi(\boldx + h \bolde_i) - \phi(\boldx)} h
\\
\Dmi \phi (\boldx) & = &
	\frac {\phi(\boldx) - \phi(\boldx - h \bolde_i)} h
\\
\Dci \phi (\boldx) & = &
	\frac {\phi(\boldx + h \bolde_i) - \phi(\boldx - h \bolde_i)} {2h}
\\
\label {eq:differenceOpsLast}
\boldD0 & = &
	(D^0_1, D^0_2, D^0_3)
\end{eqnarray}
where $ i = 1, 2, 3 $, and $\bolde_1$, $\bolde_2$, $\bolde_3$ form
an orthonormal basis of $\Rthree$.

In step 3 we use the already known $\boldu^n$ and  $\boldF^n$ to
compute  $\boldu^{n + 1}$ and $p^{n+1}$ by solving the following
linear system of equations:

\begin {eqnarray}
\label {eq:discreteNS1}
\rho \left( \frac {\boldu^{n + 1} - \boldu^n} {\Delta t}
        + \sum_{k = 1}^3 u_k^n D_k^\pm\boldu^n \right)
        & = &
\mbox{} 
        - {\bf D}^0 p^{n + 1}
        + \mu \sum_{k = 1}^3 D^+_k D^-_k \boldu^{n + 1}
        + \boldF^n
\\
\label {eq:discreteNS2}
{\bf D}^0 \cdot \boldu^{n + 1} & = & 0
\end {eqnarray}
Here 
% in equation (\ref{eq:discreteNS1}) 
% the notation 
$u_k^n D_k^\pm$
stands for upwind differencing:
$$
u_k^n D_k^\pm = \left\{
\begin {array} {ll}
u_k^n D_k^- 
& 
u_k^n > 0
\\
u_k^n D_k^+ 
&
u_k^n < 0
\end {array}
\right.
$$
Equations (\ref{eq:discreteNS1}) and (\ref{eq:discreteNS2})
are linear constant coefficient difference equations 
and, therefore, can be solved efficiently with the use of the Fast
Fourier Transform algorithm.

We now turn to the discretization of equations
(\ref {eq:fluidF}), (\ref {eq:fluidU}).
Let us assume, for simplicity, that $\Omega \subset \Rtwo$
is a rectangular domain over which all of the quantities related
to the shell are defined. We will assume that this domain is
discretized with mesh widths $\Delta q_1$, $\Delta q_2$
and the computational lattice for $\Omega$ is the set
\def\boldQ{{\bf Q}}
$$
\boldQ = \left\{(i_1\Delta q_1, i_2 \Delta q_2) \ | 
\ i_1 = 1 \ldots n_1,
\ \ 
i_2 = 1 \ldots n_2 \right\}.
$$

In step 2 the force $\boldF^n$ is computed using the following equation.
\begin {eqnarray}
\label {eq:discreteF}
\boldF^n(\boldx) = \sum_{\boldq \in \boldQ}
	\boldf^n(\boldq)\delta_h(\boldx 
	- \boldX^n(\boldq)) \Delta\boldq
\end {eqnarray}
where 
% $\boldq = (q_1, q_2)$,
$ \Delta\boldq = \Delta q_1 \Delta q_2 $
and 
$\delta_h$ is a smoothed approximation to the Dirac delta function
on $\Rthree$ described below. 
% Here the summation is over the values
% $ \boldq = (i \Delta q_1, j \Delta q_2) $ 
% where $i$ and $j$ are integers.

Similarly, in step 4 updating the position of the immersed material 
$\boldX^{n + 1}$ is done using the equation
\begin {eqnarray}
\label {eq:discreteX}
\boldX^{n + 1}(\boldq) = \boldX^n(\boldq) + \Delta t 
\sum_\boldx \boldu^{n + 1}(\boldx) \delta_h(\boldx - \boldX^n(\boldq)) h^3
\, ,
\end {eqnarray}
where the summation is over the lattice
$ \boldx = (hi, hj, hk) $, where $i$, $j$ and $k$ are integers.

The function $\delta_h$ which is used in 
(\ref {eq:discreteF}) and (\ref {eq:discreteX}),
is defined as follows:
$$
\delta_h(\boldx) = 
	h^{-3}\phi(\frac {x_1} h) \phi(\frac {x_2} h) \phi(\frac {x_3} h)
\, ,
$$
where
$$
\phi(r) = \left\{
\begin {array} {cl}
\frac 1 8 (3 - 2 |r| + \sqrt{1 + 4|r| - 4 r^2}) 
&
|r| \le 1
\\
\frac 1 2 - \phi(2 - |r|)
&
1 \le |r| \le 2
\\
0 & 2 \le |r|
\end {array}
\right.
$$
For an explanation of the construction of $\delta_h$ see 
\cite {Peskin1994}.

\subsection{Implementation and Optimization of Immersed Boundary Computations}
% \subsection{Optimization and Performance}

Our main simulation code implements the first-order
immersed boundary algorithm described above.
It is written in C and 
% The immersed boundary code 
has been optimized to run 
on several platforms: 
the Silicon Graphics Origin 2000 parallel architecture, the Cray T90
vector-parallel machine and the HP V-class and Superdome systems. 
There is also a version for a PC running Windows, suitable for
testing small models. 
The simulation code takes as input a set of files which
contain the description of the geometry and the material properties
of the system. 
These files are created before the start of the simulation by
the model generation software.
% already described.
Once the model has been read, 
the simulation enters a computation loop over
the required number of time steps,
generating output files containing information about the state of the
simulated model.

The complexity of an immersed boundary computation
is determined by
the sizes of the fluid and immersed boundary grids 
and by the size of the time step.
To prevent fluid leakage 
the mesh width of the material grids is taken to equal
approximately half the fluid grid mesh width.
(For more on volume conservation in immersed boundary calculations
see \cite{Peskin1993}).
The computation of the material forces is relatively inexpensive in
time. 
The bulk of the computation is related
to the fluid equations and to the fluid-immersed material interaction.
The solution of the discretized fluid equations 
(\ref{eq:discreteNS1}), (\ref{eq:discreteNS2})
uses the Fast Fourier Transform algorithm.
% The fluid solver relies primarily on efficient parallel 
% Fast Fourier Transform routines.
The other two demanding parts of each time step required development
of a specialized algorithm that carefully partitions the fluid and the
material grids into portions that are distributed amongst the
available processors. 
This is done in such a way as to ensure that 
no two processors operate simultaneously 
on the same portion of the data.
Numerical consistency conditions require
reducing the time step when
the space mesh width is decreased.
A convergence study of the algorithm shows that
the change in the time step is proportional to the change in mesh width.
For example, decreasing the mesh width by a factor of 2
necessitates decreasing the time step
by approximately a factor of 2 as well.
Thus rescaling a model with a $128^3$-point fluid grid to
a $256^3$-point model
requires approximately 8 times as much memory
and 16 times as much CPU time
\cite{G}.
% (see \cite{G}, chapter 4).

For large scale simulations, such as that of the full cochlea model,
extensive optimization of the immersed boundary code was necessary 
so as to reduce the elapsed time to a length that allowed useful
experiments to be completed in a manageble time. 
The code is instrumented with calls to system timing routines,
this information proving invaluable during porting 
and tuning of the software on 
different architectures. 
The largest and most recent immersed boundary 
numerical experiments have been carried out on the 32 processor 
HP Superdome installed 
at the Center for Advanced Computing Research (CACR) at Caltech.
% \footnote{
The HP 9000 Superdome at CACR contains 32 RISC processors arranged 
in a cell-based hierarchical 
crossbar architecture, 
with each cell consisting of 4 cpus with 8Gb of memory
and an I/O sub-system. 
This architecture supports the shared memory programming model
and the code efficiency was achieved primarily through 
the use of OpenMP parallelization directives.
% }
%
We used several software tools such as 
the HP CXperf and the KAI Guideview to
examine the parallel efficiency of every section in the code 
and to identify 
data cache and TLB misses. 
Figure 
% \ref{fig:cpu} 
1
shows the excellent scaling behavior of the
simulation as a function of the number of allocated processors.
Following this work, the ``wall-clock'' time per step of
the simulation is approximately 1.3 seconds 
when running on all processors of the HP Superdome.

\begin{figure}[h]
\begin{center}
  \includegraphics[width=12cm] {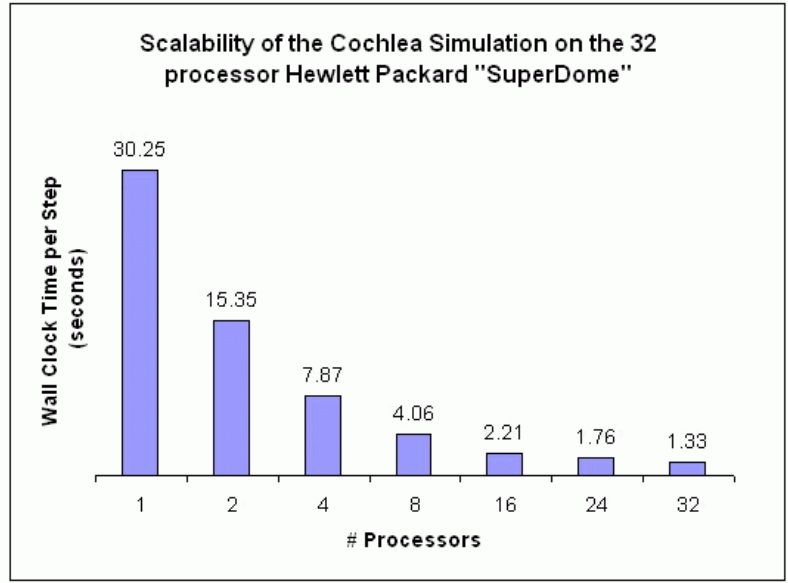}
\caption{
Execution wall-clock time per time step vs. number of processors.
}
\end{center}
\label{fig:cpu}
\end{figure}

% Each time step of the immersed boundary 
% computation consists of computing the forces in the immersed material,
% spreading these forces from the computational grid of the material
% to the computational grid of the fluid, solving the discretized
% Navier-Stokes equations and finally advancing the position of the immersed 
% material using the newly obtained fluid velocity.

%============================================================================

\section{The Cochlea Model}
\label{sec:model}

% \subsection{Anatomy and Visualization}

The size of the human cochlea is about
1\,cm $\times$ 1\,cm $\times$ 1\,cm,
while the human basilar membrane is approximately 
3.5\,cm long, is very thin and very narrow: 
its width grows from about 150\,\micron\ near the 
stapes to approximately 560\,\micron\ near the helicotrema.
Since the mesh width of the basilar membrane should approximately equal
half the mesh width of the fluid grid (see section \ref {sec:ibm}),
a fluid grid of at least $256^3$ points is necessary
to adequately resolve the width of the basilar membrane.
The geometric model of the cochlear anatomy is based on measurements
that include the position of the center line of the basilar membrane, 
its width and the
cross-sectional area of the scalae.
There are six surfaces in this model: the basilar membrane,
the spiral bony shelf, the tubular walls of the scala vestibuli and
the scala tympani 
and the semi-elliptical walls sealing the cochlear canal
and containing the oval and the round window membranes
(see Figure 
% \ref{fig:CochleaGrids}
2).
\begin{figure}[h]
\label{fig:CochleaGrids}
\begin{center}
  \includegraphics[width=14cm] {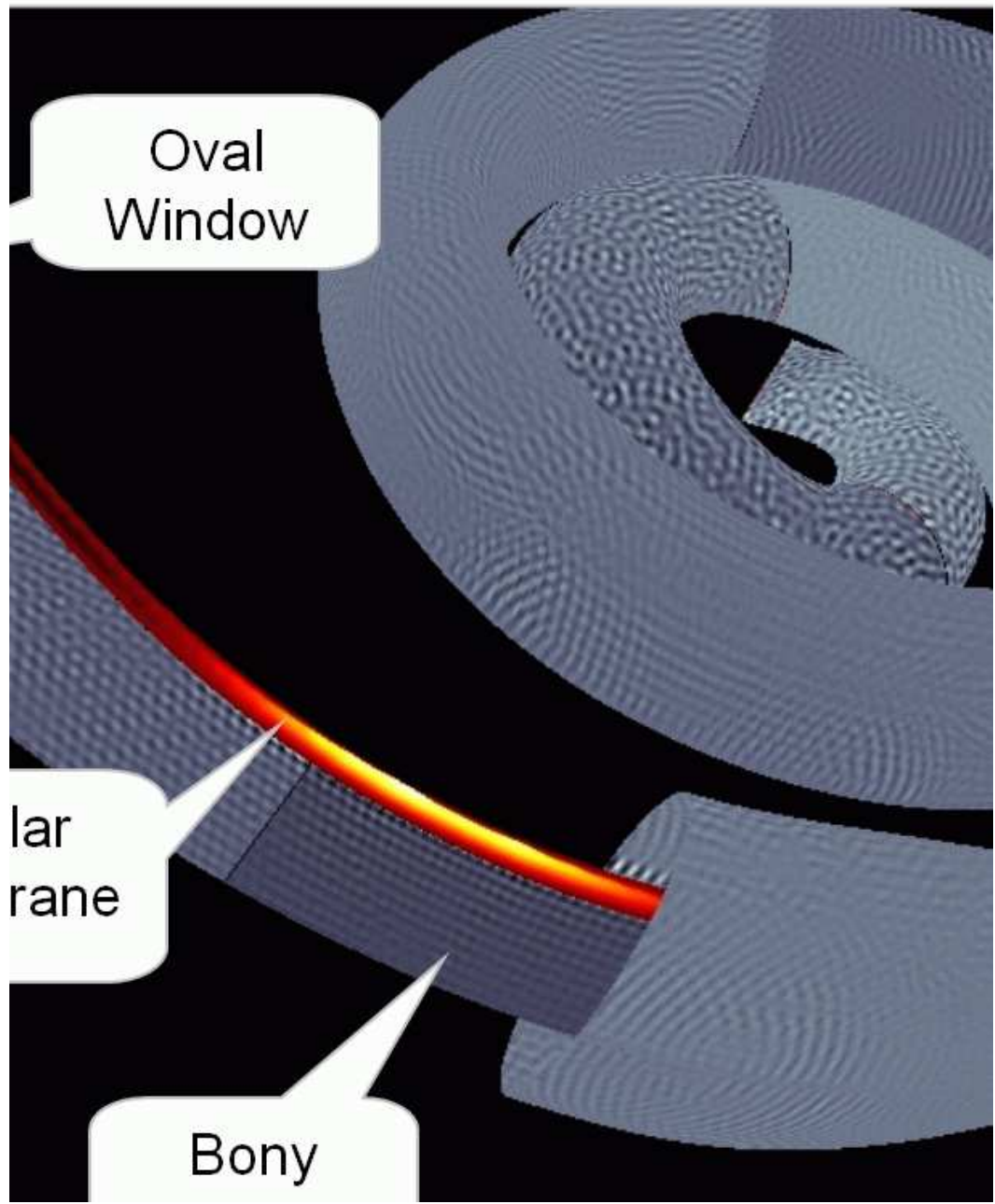}
\caption{A rendering of the geometric model of the cochlea with
several parts of the outer shell
removed in order to expose the cochlear partition consisting
of the narrow basilar membrane and the bony shelf.
The round window is located directly below the oval window
and in this picture it is partially obscured by the cochlear partition.
}
\end{center}
\end{figure}
The basilar membrane is modeled as an elastic shell following the
prototype tested earlier (see \cite{G}).
The oval window and the round window membranes are are also modeled
as elastic shells, but unlike in the case of the basilar membrane,
the compliance of each of these shells is constant throughout the
shell.
% by computational grids, whose points are interconnected by elastic
% springs.
% Geometrically, both windows are the same;
The windows are geometrically identical:
each is modeled by a rectangular grid, all of whose points outside
a given circle are fixed.
Hence each such grid represents a circular plate within a
rectangular piece of a bone.

The model building programs generate 
the approximate cochlear geometry,
the material properties of various surfaces
and the set of parameters describing the desired simulation tests.
The six surfaces of immersed material in the cochlea model 
are partitioned into
31 computational grids comprising approximately 750,000 points
in total.
In addition to the basilar membrane and the windows' grids described
above there are 28 grids modeling the bony surfaces.
The partitioning of the bony surfaces into 28 rectangular 
grids was necessary to minimize the total number of material points
while maintaining an approximately uniform distance between these
points.

% for small amplitude input 
The passive cochlea is essentially
a linear mechanical filter (see section \ref{sec:cochlea})
and the displacements occuring within the cochlea 
are too small to be seen
without magnification.
% At the hearing threshold the displacement of the basilar membrane
% are measured by a fraction of a nano-meter.
In our numerical experiments the displacements of the basilar membrane
are measured on a nano-meter scale.
The immersed boundary method is
particularly suitable for such 
a simulation since it possesses a subgrid resolution
% for example, the  calculated displacements of the \BM\ are
% on the order of nanometers,
% % with an estimated error of less than ???,
% while 
(the material grid mesh width of our cochlea model is about 20\micron).
Our numerical algorithm however uses the fully non-linear system
of Navier-Stokes equations rather than 
the linearized system
because the computational cost of
solving the linearized equations in the immersed boundary method
would be essentially the same as the cost of solving 
the fully non-linear system.

The codes which are used to analyse and visualise 
% the simulated data
the cochlear geometry and the simulated dynamics
include a C++/OpenGL tool that runs on SGI and a Java/Java3D tool 
that runs on most platforms. 
These tools take as input the vertex coordinates for
all computational grids in the cochlea model
from the result files for each time step in the simulation. 
The tools generate
a full 3D rendering of the model geometry.
Since the displacements occuring within the cochlea
% Since the displacements of the basilar membrane 
are very small
% (measured by a fraction of a nano-meter at the hearing threshold),
they are color-coded to reveal the dynamics of the system.
Our graphics tools also enable us to magnify the displacements
to make them visible on the screen.
% where each grid vertex 
% is displayed in
% a color that corresponds to its displacement from the 
% initial (starting) position.
% Optionally, the displacements along the normals of the 
% vertices from the starting  
% position can be exaggerated (scaled) by a factor of choice so as 
% to reveal
% the motion of the surface, which would otherwise be imperceptibly small in the display. 
% By rendering the surface of the basilar membrane in this way for each time step
% in the simulation, and using the option of saving each rendered frame as a
% JPEG file, we are able to afterwards build a movie by compiling together the
% series of frames that were produced. 
An essential insight into the basilar membrane dynamics is provided
by the plot of the normal displacement of its center line
(see Figures \ref{fig:bm_center1}, \ref{fig:bm_center2}).
Other Java tools display various
important characteristics of the system dynamics,
such as the response of individual points 
on the basilar membrane as a function of time.
All of the graphics tools have built-in facilities 
for generating animation.

%\section{Description of the Model Software, Performance and Validation}
% \subsection{Software Components}
% \subsubsection{Model Generation Software}
% \subsubsection{Immersed Boundary Numerical Solver}
% \subsubsection{Analysis and Visualization}

% \subsection{Validation of the Model Components}

The construction of the full cochlea model has been undertaken
in stages with the individual components tested separately 
prior to the model assembly.
The basilar membrane model is described in detail in \cite{G}.
The oval window functions as the input window of the cochlea.
We simulate the pressure of the stapes against it
by prescribing an external force vector field on the window grid.
This force is orthogonal to the surface of the window
and for a pure tone input its magnitude varies sinusoidally.
No force is prescribed on the round window.
%% % We have examined the resulting motion of the windows
%% % and a snapshot image from our Java 3D application
%% % is shown in Figure \ref {fig:DualWindows}.
%% A snapshot image of the two windows in motion during a simulation test
%% is shown in Figure \ref {fig:DualWindows}.
%% We used our Java3D application to exaggerate and color-code
%% vertical displacements.
%% 
%% \begin{figure}[h]
%% \begin{center}
  %% \includegraphics[width=12cm] {DualWindows.eps}
%% \caption{A snapshot of the oval and the round windows 
	%% in motion during a simulation test. 
	%% The Java3D application is used to exaggerate and color-code
	%% vertical displacements.}
%% \end{center}
%% \label{fig:DualWindows}
%% \end{figure}

Since many of the numerical experiments with the full cochlea model
require days of computing we have made a small modification
of the cochlear anatomy in our model.
% Several small changes were 
% made to the geometric model of the cochlear anatomy. 
In the altered model 
the cochlea is packed more tightly and the whole structure fits in a 
half-cube of size 
1cm $\times$ 1cm $\times$ 0.5cm
rather than a 1 cm$^3$ cube, making it possible to use a fluid grid of 
$256 \times 256 \times 128$ points.
This configuration requires about 1.5 Gigabytes of main memory
and the total time needed to simulate a single time step
is reduced approximately in half.
The main difference between the ``half-cube model'' 
and our original model is that in the former model the basilar membrane
if flat. I.e. its spiraling shape is completely contained within
a plane.
All of our test simulations reported in the next section have been
carried out with the half-cube model.

% Since the cochlea model we are using is comprised of several different
% component parts (the Oval Window, the Round Window, the Basilar Membrane,
% etc.) we were able to check the behavior of each component piecewise before
% assembling them into the complete model.
% 
% The behavior of the Oval Window when excited by a sinusoidal force was first
% examined. 
% We ran experiments with various input frequencies and amplitudes,
% and verified that a choice of 30ns as a time step size was optimal. 
% Too small a time step resulted in numerical instability in the solver, 
% and too
% large a time step resulted in the model results not showing the small scale
% behavior of the system correctly. 
% Figure X shows a snapshot image from our Java 3D
% application that renders the Oval Window displacement as a function of time.
% Note that the force is applied over a circular region of the Oval Window. The
% snapshot is taken after X milliseconds of simulated time.

% The principal practical difficulty in the implementation 
% is due to the fact that the
% basilar membrane is very narrow relative to the size 
% of the whole cochlea 
%

\section{Numerical Experiments}
\label{sec:results}

The good efficiency of the immersed boundary solver
allowed us to complete several large numerical experiments.
% which we describe in the remainder of this paper.
%
In this section we present preliminary results from four such
experiments.
% We present some preliminary results in this section.
% In this section we describe some of the results obtained.
We are continuing this work and
a more detailed exposition of the collected data and its analysis
will be published in our next paper.

Each of the experiments presented here consisted of applying 
a pure tone input
of a given frequency
at the stapes 
and studying the subsequent motion of the basilar membrane.
The input at the stapes was simulated by specifying an 
external periodic force vector field on the oval window grid.
A very small time step of 30 nano-seconds was chosen to guarantee
numerical stability and good detail.
%, as explained above.
The choice of the time step was made as a result of 
the convergence study of the system 
carried out in \cite{G}.
% (this is discussed in detail in \cite{G}).
%
Our initial aim in the numerical experiments was to reproduce the
qualitative characteristics of cochlear mechanics 
predicted by
asymptotic analysis and previously reported 
computational models.
Consequently, we have attempted to capture the steady state response
of the basilar membrane to pure tones and 
in each of our experiments we have run the system for the duration
equivalent to several input cycles.
For example, for a 10 kHz input tone we have run the system
for more than 30,000 time steps,
% With the simulation time step of 30 nano-seconds this is 
equivalent to 0.9 msec of total simulated time.
On a 16 processor SuperDome
this computation required more than 20 hours to complete.
Every 10th time step of the simulation the program
generated an output file containing the instantaneous position
of the computation grids.
The total amount of storage required is measured in tens of Gigabytes
of disk space.
Correspondingly, the 2 kHz experiment required five times as much
computational time and storage.
% A typical experiment was run for 50,000 time steps,
% equivalent to 0.2 milliseconds of simulated time,
% and took approximately 18 hours of dedicated computation
% on the Superdome.
%

\begin{figure}[h]
\begin{center}
  \includegraphics[width=10cm] {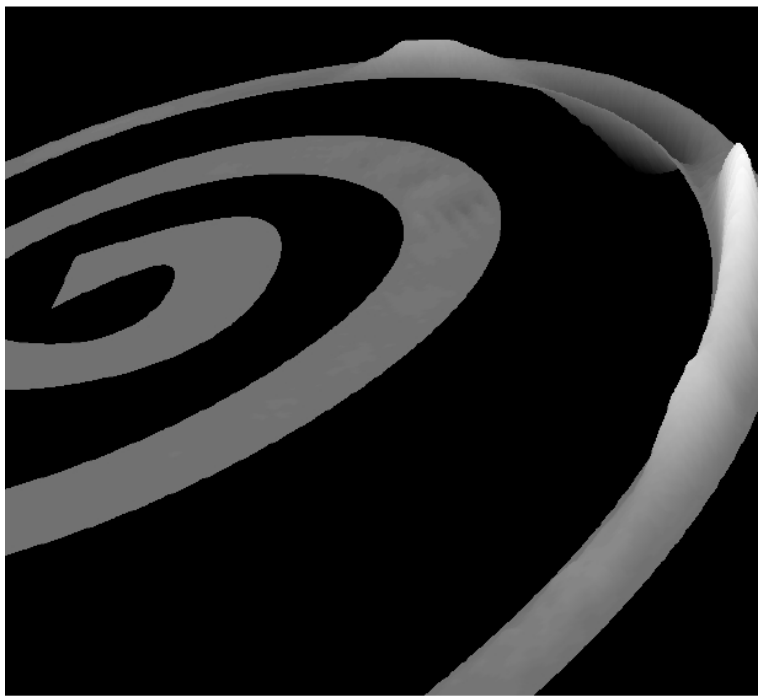}
\caption{
A close-up snapshot 
of the traveling wave propagating along the basilar membrane.
The magnitude of the basilar membrane displacement 
has been amplified in the normal direction.
}
\end{center}
\label{fig:wave}
\end{figure}

We have verified that all of our experiments have been carried in
the system's linear regime, i.e. the input force at the stapes was kept
very small, resulting in the basilar membrane displacements on
the order of nano-meters.
Indeed, increasing the force by a factor of 10 resulted in
basilar membrane displacements almost exactly ten times bigger.
Since in each experiment the system was started from rest,
we have observed initial oscillatory transient response,
which was followed by the smoothing of the traveling wave.
A close-up snapshot of the traveling wave
propagating along the
basilar membrane is shown in Figure 
% \ref{fig:wave}.
3.
The wave magnitude has been amplified in the direction of the normal to
the basilar membrane.

\begin{figure}[h]
% \setlength{\epsfxsize}{5in}
% \hspace{0.5in}
\vspace{0.2in}
% \epsffile{bm_10F_center.3500.ps}
\resizebox{!}{2.5in}{\includegraphics{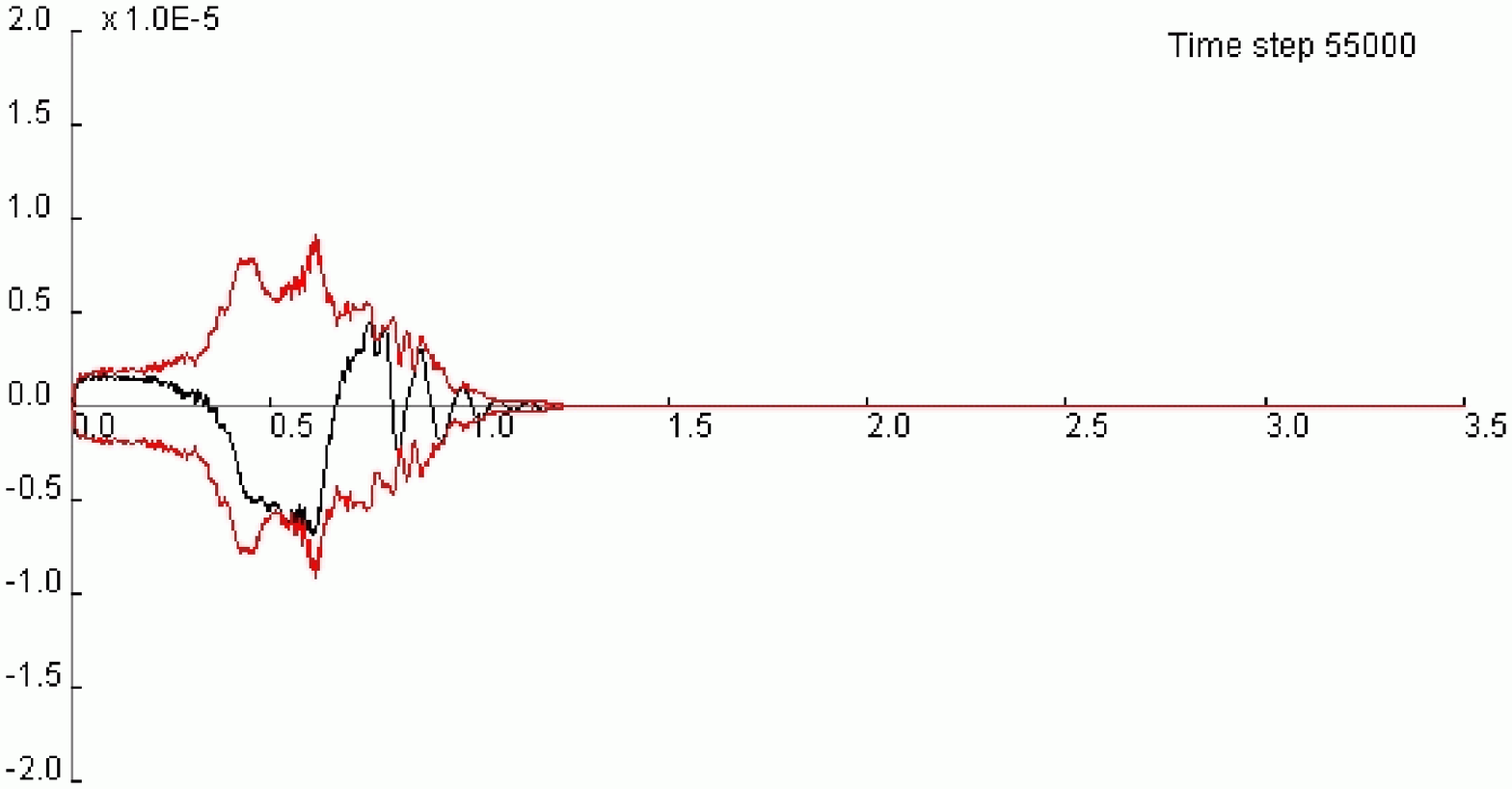}}
\resizebox{!}{2.5in}{\includegraphics{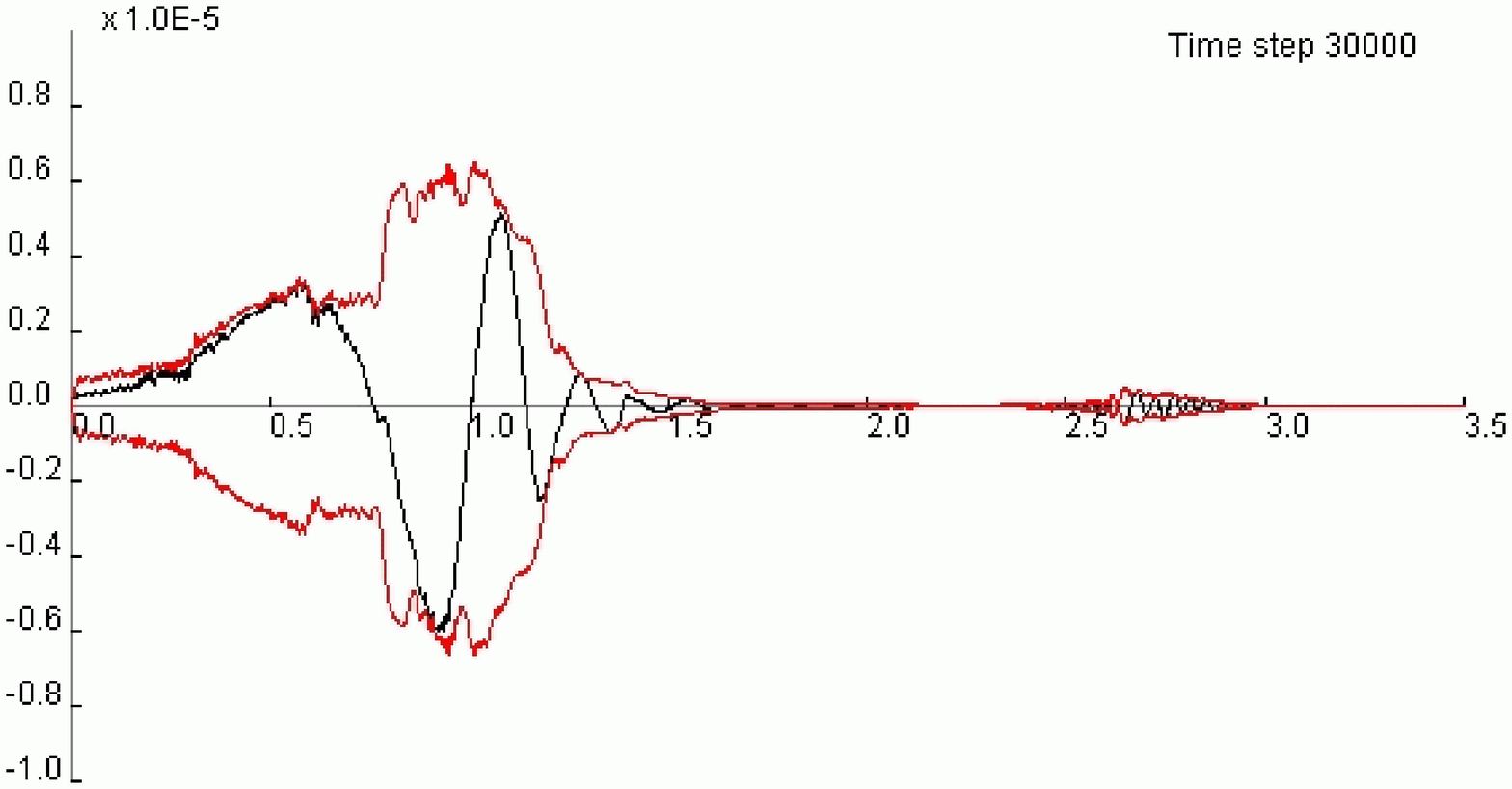}}
\caption{
A snapshot of the center line of the basilar membrane
response and wave envelope.
The top plot shows the response to the 15 kHz input sound.
The wave snapshot is taken after 55000 time steps
(1.65 msec).
The wave envelope was computed over time steps
45,000 - 55,000.
The bottom plot shows a 10 kHz experiment.
The wave is shown after 30,000 time steps and the envelope was computed
over time steps 20,000 - 30,000.
The unit of vertical scale is 10$^{-5}$cm.
}
\label{fig:bm_center1}
\end{figure}

\begin{figure}[h]
% \setlength{\epsfxsize}{5in}
% \hspace{0.5in}
\vspace{0.2in}
% \epsffile{bm_10F_center.3500.ps}
\resizebox{!}{2.5in}{\includegraphics{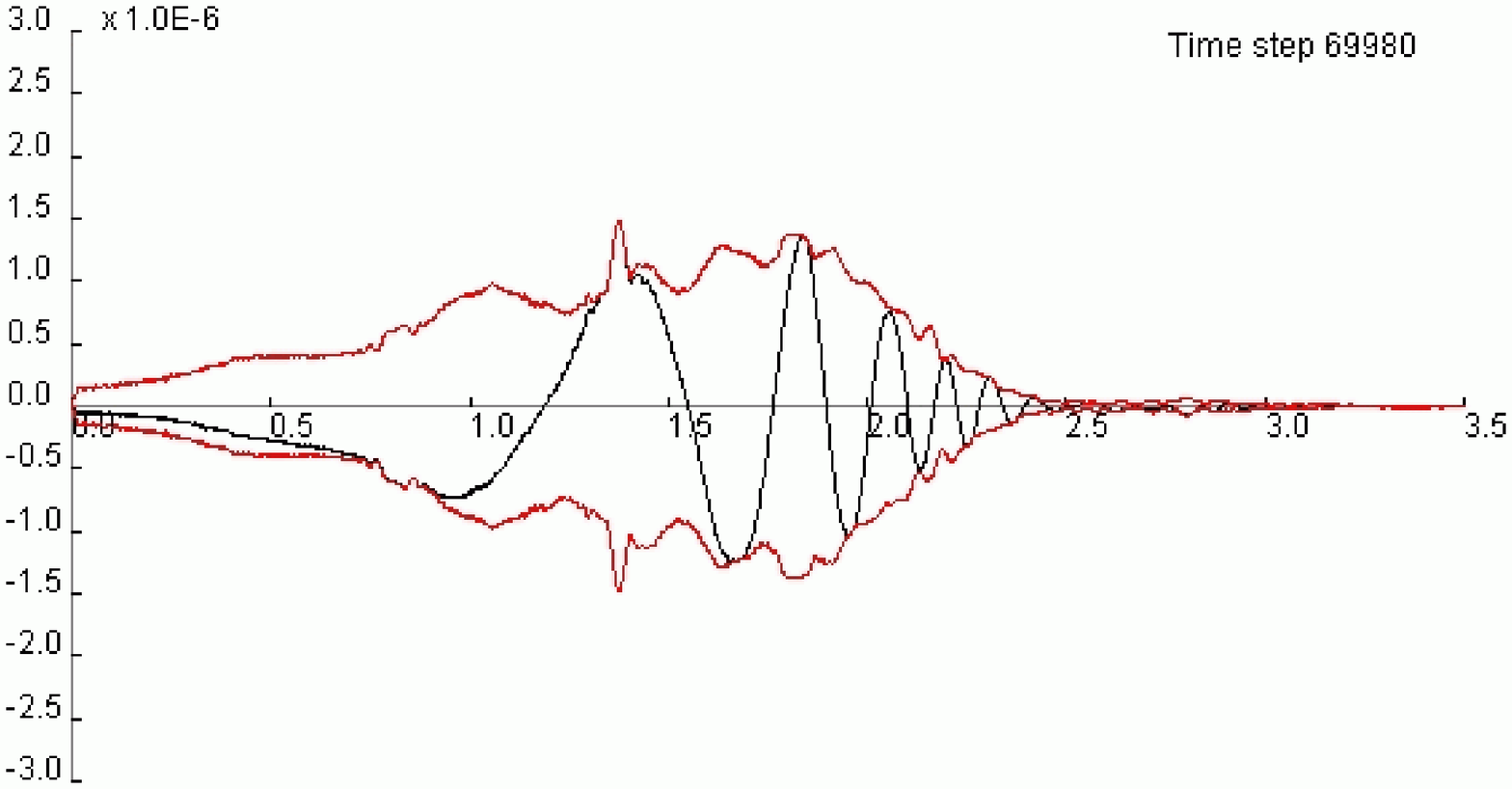}}
\resizebox{!}{2.5in}{\includegraphics{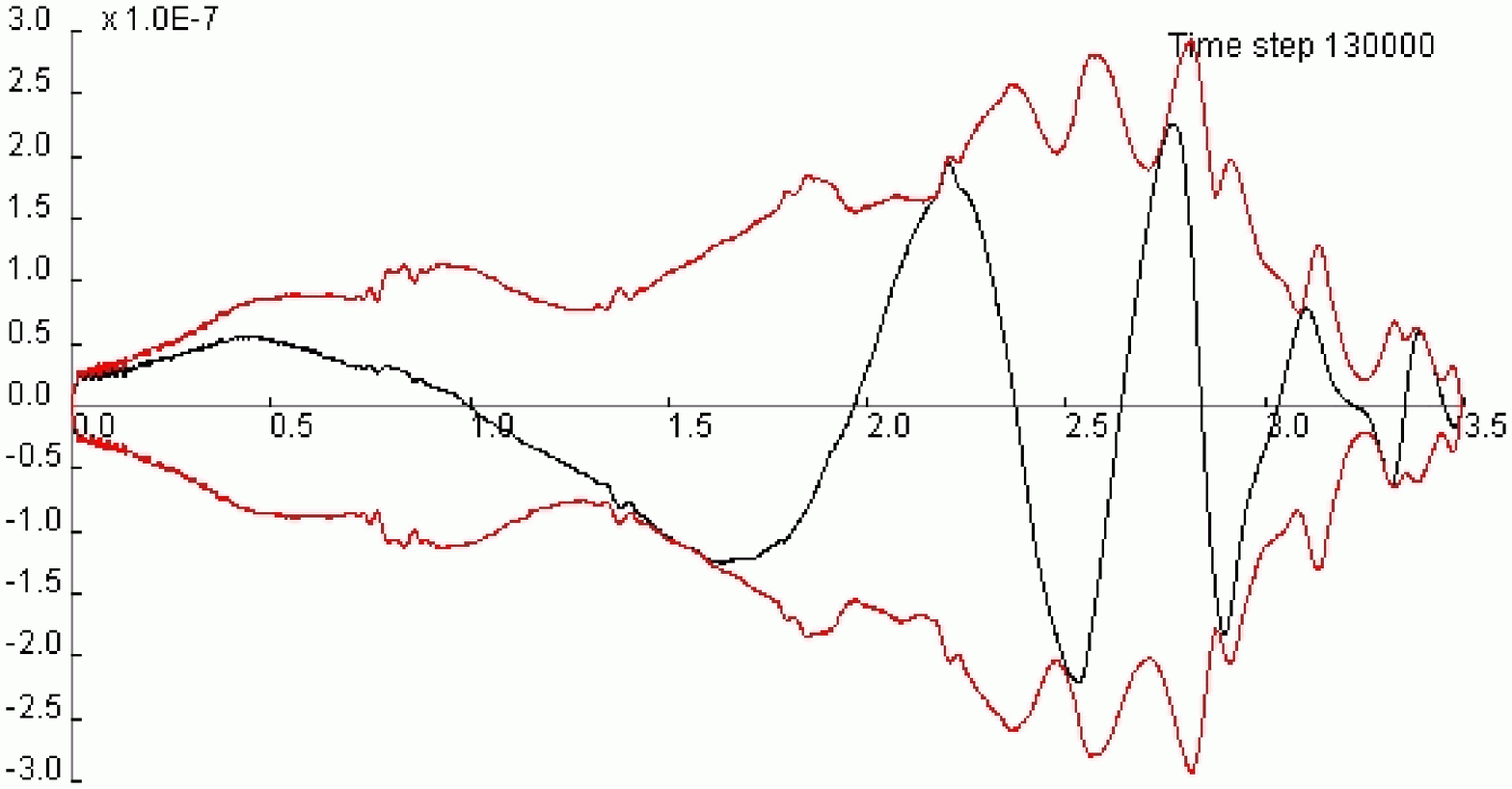}}
\caption{
A snapshot of the center line of the basilar membrane
response and wave envelope.
The top plot shows the response to the 5 kHz input sound.
The wave snapshot is taken after 70000 time steps.
The wave envelope was computed over time steps
60,000 - 70,000.
The bottom plot shows a 2 kHz experiment.
The wave is shown after 100,000 time steps and the envelope was computed
over time steps 100,000 - 130,000.
The unit of vertical scale is 10$^{-6}$cm and 10$^{-7}$cm in the top
and in the bottom plots, respectively.
}
\label{fig:bm_center2}
\end{figure}

Much information about the traveling wave is revealed 
in the plot of the centerline of the basilar membrane.
Four such plots are shown in Figure \ref{fig:bm_center1}
and Figure \ref{fig:bm_center2}.
Figure \ref{fig:bm_center1} 
shows the results of the experiments with input frequency
of 15 kHz (top plot) and 10 kHz (bottom plot),
and Figure \ref{fig:bm_center2} shows the results of the
5 kHz and the 2 kHz experiments.
In each plot an instantaneous position of the centerline
of the basilar membrane is shown,
as well as the wave envelope
computed over a range of time steps.

Our experiments reproduced the following characteristic features
of cochlear mechanics:
In each instance,
in response to a pure tone input frequency,
we have observed a traveling wave propagating from the stapes in the
direction of the helicotrema.
The amplitude of the wave is gradually increasing until it
reaches a peak at a characteristic location along 
the basilar membrane depending on the
input frequency.
The speed of the wave is sharply reduced
as it propagates further along the basilar membrane.
% This is to be expected in view of the increasing compliance
% of the basilar membrane.
The higher the input frequency, the closer the peak of the wave
is to the stapes.
Furthermore, 
after reaching the peak the wave drops off
sharply, essentially shutting down.
Since no active mechanism has been incorporated into the present model
yet,
the observed traveling wave is not sharply focused.

It is interesting to note that while the computed traveling wave 
is smooth, its envelope is generally not smooth.
% We have observed this in every one of the four numerical experiments.
The 10 kHz experiment shown in Figure \ref{fig:bm_center1}
is a particularly striking example,
with the envelope turning 90$^\circ$
sharply in the vertical direction approximately at the
0.75$~cm$ location along the basilar membrane.
We believe that further testing is necessary 
and more data must be collected from
numerical experiments before further conclusions about the shape of the
traveling wave envelope can be drawn.

% Excitation using a pure tune sinusoidal input was employed to
% pinpoint the system response at the frequency used. 
% The behavior of
% the basilar membrane is observed to be quite different in this case:
% a sharp oscillation at one point along the basilar membrane is produced, in
% contrast to the "click" response, where a wave packet of decreasing amplitude
% and velocity travels along the membrane.

% The speed of the wave packet is sharply reduced
% as it propagates along the basilar membrane.
% This is to be expected in view of the increasing compliance
% of the basilar membrane.
% The results of these experiments will be 
% compared with Peskin's asymptotics
% \cite{Leveque1988}
% and
% Karl Grosh's hybrid analytical-computational model
% \cite{Parthasarati2000},
% and reported in future publications.
%
The interested reader is invited to view several animations 
of our results by visiting our web site\cite{Julian_website}.

\section{Summary and Conclusions}
\label{sec:future}

We have constructed a comprehensive three-dimensional computational
model of the passive cochlea using the immersed boundary method.
Extensive optimization and parallelization made it 
possible to complete several large scale numerical simulations
on a 32-processor shared memory HP Superdome computer.
Together with the previous demonstration of the traveling wave paradox
(see \cite{G}),
the pure tone experiments reported in this paper
capture the most important properties of
the cochlear macro-mechanics.

We would like to note that our results are of somewhat 
preliminary nature.
We are continuing to test the model and plan to complete
more numerical experiments in order to compare our simulation
results with the available experimental and modeling data.
We believe our results demonstrate the promise of large scale
computational modeling approach to the study of cochlear mechanics.

The cochlea is a very small and delicate organ
and it is very difficult to study experimentally.
Many of the important questions of cochlear mechanics 
are mathematically very complicated,
but they can be studied using numerical simulation.
Two such general questions concern the effects of the geometry 
and of the elastic properties of various cochlear components
on the dynamics.
It is easy to test cochlear models with different geometries
such as a straightened out model.
Testing different elastic models for the cochlear components
is somewhat more involved since 
it generally requires testing each such model separately within
the immersed boundary framework.
It is nevertheless a feasible project, which is straightforward 
for the expert.

Going beyound the questions of passive cochlear macro-mechanics 
we would like to incorporate an active mechanism into our model.
The function of the live cochlea undoubtedly depends on the
combination of passive cochlear mechanics with an active mechanism.
A comprehensive model of the passive cochlea is therefore a
necessary first step towards modeling the active live cochlea.
It is interesting to note that refining the mesh width of the 
computational grids by a factor of 2 yields a material mesh width
of approximately 10\micron.
This mesh width is small enough to allow much of the structure of the
organ of Corti to be incorporated into the model.
Such a refined model is not yet feasible since
it would require approximately 16 times more computing
power than our present model.
Indeed, the large scale of computations is presently
the biggest obstacle 
to progress in cochlear modeling,
but this is, no doubt, a temporary obstacle.
Continuing progress in hardware and software will make
a construction of even the refined cochlear model possible soon.

% We intend to continue testing the model and compare the results
% with the available experimental and modeling data.
% For example, we would like to construct a cochlear map for this model,
% i. e. to determine the characteristic frequency of each location
% on the basilar membrane.
% Subsequently we will extend the present model to capture 
% the micro-mechanics of the organ of Corti in order to study
% the microscopic mechanism of amplification in the cochlea.
% We anticipate such a model to be even more computationally demanding
% than the present model.

\section{Acknowledgements}
The first author would like to thank Charlie Peskin 
for suggesting the project
and for his guidance through its initial phase,
and Karl Grosh for many helpful discussions.
We would like to thank Mahesh Rajan for his help with optimizing the
numerical solver on the HP machines, and Sarah Emery Bunn for her 
assistance in the preparation of the manuscript.

\raggedright
\bibliographystyle{plain}
\bibliography{cochlea}

\end{document}